\documentclass[11pt]{article}

\usepackage{amsmath}
\usepackage{amsfonts}
\usepackage{amssymb}
\newtheorem{theorem}{Theorem}
\newtheorem{lem}[theorem]{Lemma}
\newtheorem{cor}[theorem]{Corollary}
\newtheorem{prop}[theorem]{Proposition}
\newtheorem{conj}[theorem]{Conjecture}

\def\Z{\mathbb Z}

\begin{document}
\author{K\'aroly J. B\"or\"oczky\footnote{Supported by
OTKA grants 068398 and 75016, and by the EU Marie Curie FP7 IEF grant 
GEOSUMSETS},
P\'eter P. P\'alfy\footnote{Supported by
OTKA grant NK72523}, Oriol Serra}

\title{On the cardinality of sumsets  in torsion-free groups}

\date{}

\maketitle

\begin{abstract}
Let $A, B$ be finite subsets of a torsion-free group $G$. We prove
that for every positive integer $k$ there is
a $c(k)$ such that if
$|B|\ge c(k)$ then  the inequality $|AB|\ge |A|+|B|+k$ holds unless
a left translate of $A$ is contained in a cyclic subgroup. We obtain
$c(k)<c_0k^{6}$ for arbitrary
torsion-free groups, and $c(k)<c_0k^{3}$
for groups with the unique product property, where $c_0$
is an absolute constant. We give examples to
show that $c(k)$ is at least quadratic in $k$.
\end{abstract}

\section{Introduction}

Let $G$ be a torsion-free group written multiplicatively,
and let $|\cdot|$ denote the cardinality
of a finite set.  A basic
problem in Additive Combinatorics is to estimate the cardinality of
$AB=\{ ab:\; a\in A, b\in B\}$ of two finite sets $A,B\subset G$ in
terms of $|A|$ and $|B|$.
A basic notion is a progression with ratio $r\neq 1$
and length $n$, which is a set of the form $\{a,ar,\dots,ar^{n-1}\}$
where $a$ and $r$ commute.

Let us review some related results if $G$ is abelian.
In this case we have the simple inequality
$$
|AB|\ge |A|+|B|-1,
$$
with equality if and only if
 $A$ and $B$ are progressions with common ratio.
Following Ruzsa \cite{Ruz09},
we call the minimal rank of
a subgroup whose some coset contains  $A$
the dimension of $A$.
According to Freiman \cite{Fre73}, if the dimension of $A$
 is $d$, then
\begin{equation}
\label{Freimanmulti}
|A^2|\ge (d+1)|A|-{d+1 \choose 2}.
\end{equation}
This estimate is optimal.
It follows that if $|A^2|\leq 3|A|-4$, then
$A$ is contained in some coset of a cyclic group.
Actually, even a progression of length
$2|A|-3$ contains $A$ according to the
 $(3k-4)$-theorem of  Freiman \cite{Fre73}.
More precise structural information on $A$ is available
 if $|A^2|= 2|A|+n$
for $0\leq n\leq |A|-4$ by  Freiman \cite{Fre09b},
for example, $A$ is contained in a progression of length $|A|+n+1$.

The inequality (\ref{Freimanmulti}) was generalized to
a pair of sets by Ruzsa \cite{Ruz94}
who proved that if  $|A|\ge |B|$, and the dimension
of  $AB$ is $d$, then
\begin{equation}
\label{Ruzsa}
|AB|\ge |A|+d|B|-{d+1\choose 2}.
\end{equation}
 By requiring additionally that the smaller set $B$ is
$d$--dimensional, Gardner and Gronchi \cite{GaG01} proved a
discrete version of the Brunn--Minkowski inequality which shows that
\begin{equation}
\label{GG}
|AB|\ge |A|+(d-1)|B|+(|A|-d)^{(d-1)/d}(|B|-d)^{1/d}-{d\choose 2}.
\end{equation}
Additional lower bounds with stronger geometric requirements on the
sets $A$ and $B$ have been also obtained by Matolcsi and Ruzsa
\cite{MaR} and Green and Tao \cite{GrT06}.

In the non-abelian case the situation is much less
understood. Kempermann \cite{Kem56}  implies in the case
of  any torsion-free group $G$ that
\begin{equation}
\label{eq:kemp}
|AB|\ge |A|+|B|-1.
\end{equation}
 Brailovsky and Freiman \cite{BrF90} characterized the extremal sets in the
inequality \eqref{eq:kemp} by showing that, if $\min\{ |A|,|B|\}\ge
2$ then, up to appropriate left and right translations, both $A$
and $B$ are progressions with common ratio.
In particular, $A$ and $B$ lie in a left and
a right coset, respectively,
of a cyclic subgroup.

The analogy with the abelian case was extended in
Hamidoune, Llad\'o and Serra \cite{HLS98} to
the inequality
\begin{equation}
\label{HLS}
|AB|\ge |A|+|B|+1,
\end{equation}
if $|B|\geq 4$, and $A$ is not contained
in some left coset of a cyclic subgroup.

These known facts are connected with the following conjecture of
Freiman (personal communication),
extending the $(3k-4)$-theorem above.

\begin{conj}
\label{conj:freiman}
 Let $A$ be a
 finite  subset of a torsion-free group with $|A|\ge 4$. If
 $$
 |A^2|\le 3|A|-4
 $$
 then $A$ is covered by a progression of length at most $2|A|-3$.
 \end{conj}

  By using the so--called isoperimetric method, see Hamidoune
 \cite{Ham96},  \cite{Ham08} or  \cite{Ham09},
 we obtain the following results:

 \begin{theorem}
\label{thm:main}
 For any integer $k\geq 1$ there exists
a $c(k)$ such that
  the following holds.
If $G$ is a torsio-free group,
 $A\subset G$ is
  not contained in a left coset of any cyclic  subgroup,
$B\subset G$ has more than $c(k)$ elements,
then
  $$
  |AB| > |A|+|B|+k.
  $$
  \end{theorem}
{\bf Remark }  Our current methods yield $c(k)\le 32 (k+3)^6$.\\

  Note that, in Theorem~\ref{thm:main},
the assumption on $A$   not being contained in a left coset
  of a cyclic group  is crucial. For example,
  if $A$ is a progression of length at most $k+2$ with ratio $r\neq
  1$, 
and $B$ is
 the union of two $r$-progressions of arbitrary length, then 
$ |AB| \leq |A|+|B|+k$.

The value of the lower bound $c(k)$ can be
  improved for unique product groups. Recall that a
   group $G$ has the unique product property if, for every pair of
  finite sets $A, B\subset G$, there is an element $g\in AB$ which
  can be uniquely expressed as a product of an element of $A$ and an
  element of $B$. In this case, $G$ is torsion-free. We note that
every right linearly orderable
group has the unique product property, and any residually finite 
word hyperbolic group
has a finite index unique product subgroup, according to 
T. Delzant \cite{Del97}. On the other hand,
it was first shown by Rips and Segev \cite{RiS87} that not all 
torsion-free groups have the
  unique product property, and S.D. Promislow \cite{Pro88}
even provided an explicit construction for such an example.
In addition, unique product groups are discussed in
A. Strojnowski \cite{Str80}, S.M. Hair \cite{Hai03},
and W. Carter \cite{Car07}.
For unique product groups,   the bound
on $c(k)$ in Theorem~\ref{thm:main} can be
  reduced to a cubic polynomial on $k$; namely,
Lemma~\ref{unique} yields that
\begin{equation}
c(k)\leq 4(2k+3)^3 \mbox{ \ if $G$ is a unique product group.}
\end{equation}

We note that it can be deduced
with the help of (\ref{GG}) that in abelian torsion-free groups, 
the optimal order of $c(k)$ is quadratic.
The $c(k)$ in Theorem~\ref{thm:main} is at least of quadratic order
also for non-abelian unique product groups,
as the following example shows.

We consider the Klein bottle group $G_0=\langle u, v | u^{-1}vu=v^{-1}\rangle$,
and hence $vu=uv^{-1}$ and $v^{-1}u=uv$. Since $\langle v\rangle$ is a 
normal subgroup
with factor isomorphic to $\Z$, $G_0$ is a non-abelian unique product group.
Let $A = \{1,u,v\}$ and  $B = \{ u^i v^j : i,j=0,1,...,m-1 \}$
for $m\geq 1$.
Then $|B|=m^2$ and
$AB =\cup_{i,j=0,1,...,m-1}\{ u^i v^j, u^{i+1} v^j, u^i v^{j+(-1)^i} \}$, thus
$$
|AB| = m^2 + 2m=|A|+|B|+2|B|^{\frac12}-3.
$$

In line with Conjecture~\ref{conj:freiman}, we conjecture that, 
if $A=B$, then  the lower
bound on $|B|$ in Theorem~\ref{thm:main} can be replaced by a bound
linear in $k$.
We construct an example in the group $G_0$ above
to indicate, what to expect
in Theorem~\ref{thm:main} in this case.

 For $m\geq 1$, let $A=P\cup vuQ$ where $P=\{ u^i  : i=0,1,...,2m \}$ and
$Q=\{ u^{2i}  : i=0,1,...,m-1 \}$, and hence  $|A|=3m+1$. Since $v$
commutes with $u^2$, we have $(vuQ)(vuQ)=(vuv)uQ^2=u^2Q^2\subset
P^2$.
Moreover, denoting by $P_0$ and $P_1$ the set of even
and odd powers of $u$ in $P$, respectively, we have $PvuQ=vuP_0Q\cup
v^{-1}uP_1Q\subset vuQP\cup v^{-1}uP_1Q$. It follows that
$$
|A^2|=|P^2|+|QP|+|P_1Q|=10m-1=\mbox{$\frac{10}{3}|A|-\frac{13}{3}$}.
$$
We note that the above example seems to match  a conjecture of
Freiman \cite{Fre02}, which would yield that $A$ is the union of
two progressions provided that $|A^2|< \frac{10}3\,|A|-5$.

In the direction of Conjecture~\ref{conj:freiman} for a
torsion-free group $G$, our results yield the following.

 \begin{cor}
\label{A=B} If $A$ is a subset of a torsion-free group with
$|A|\geq 6^6$, and $|A^2|=2|A|+n$ for $0\leq n\leq
2^{-5/6}|A|^{1/6}-3$, then $A$ is contained in a progression of
length $|A|+n+1$.
\end{cor}
{\bf Remark }  In unique product groups,
the conditions are
$|A|\geq 6^3$ and
 $0\leq n\leq 2^{-5/3}|A|^{1/3}-\frac32$.\\

If $A$ is a finite subset of a torsion-free
group $G$, then Corollary~\ref{A=B}
provides strong structural information
when $|A^2|$ is very close to $2|A|$.
This has been made possible in part by the
known structural properties in abelian groups.
Now if $G$ is abelian and $|A^2|<K|A|$ for some $K> 3$
then still strong structural properties have
been established by Freiman \cite{Fre73}
using multidimensional progressions,
see the monograph
of Tao and Vu \cite {TaV06}
or the survey by Ruzsa \cite{Ruz09}
for recent developments. But if $G$
is any torsion-free group and $K\geq \frac{10}3$,
then $A$ may not be contained in an abelian
subgroup. Actually, it is still not
completely understood, what to
expect, in spite the results about some specific groups 
(see Breuillard and Green \cite{BeGI},
or Tao's blog \cite{Tao09}).

\section{Atoms and fragments}
\label{sec:atoms}

For this section, we fix a torsion-free group $G$.

For $n\geq 1$ and a finite non-empty set $C\subset G$, the $n$-th
isoperimetric number of $C$ is defined to be
$$
\kappa_n(C)=\min\{|XC|-|X|:\,X\subset G\mbox{ and }|X|\geq n\}.
$$
A finite set $V\subset G$ is an $n$-fragment for $C$, if $|V|\geq n$
and $|VC|-|V|=\kappa_n(C)$. In addition an $n$-fragment of minimal
cardinality is an $n$-atom for $C$.

Naturally, if $U$ is an $n$-atom for $C$, then $xU$ is also an
$n$-atom for $Cy$ for any $x,y\in G$. In what follows, we present
simple statements about atoms. For the sake of completeness, we
verify even the known ones, except for the following crucial
property of atoms, due to Hamidoune \cite{Ham96}: If $U$ is an
$n$--atom  and $F$ is an $n$--fragment for a finite nonempty
subset $C\subset G$, then
$$
\mbox{either $U\subset F$ or $|U\cap F|\le n-1$.}
$$
This property has the following  useful consequence.

\begin{cor}
\label{atomleft}
For a torsion-free group $G$ and $n\geq 1$, if $U$ is an $n$-atom
for $C\subset G$ and $g\in G\backslash 1$, then $|U\cap gU|\leq n-1$.
\end{cor}

For right translations  we have a weaker result.

\begin{lem}
\label{atomright}
For a torsion-free group $G$ and $n\geq 2$, if $U$ is an $n$-atom
for $C\subset G$ and $g\in G\backslash 1$, then
$$
|U\cap Ug|\leq \frac{n-2}{n-1}\,|U|+\frac{1}{n-1}
\le \frac{n-1}n\,|U|.
$$
\end{lem}
{\bf Remark } In particular, if $n=2$, then $|U\cap Ug|\leq 1$.\\
{\it Proof: } Let us partition $U$ into the maximal left
$g$-progressions $U_1,\ldots,U_m$, where
$U_i=\{h_i,h_ig,\ldots,h_ig^{\alpha_i}\}$, $i=1,\ldots,m$. In
particular, $U_i\cap U_jg=\emptyset$ for $i\neq j$. We may assume
that $|U_1|\geq |U_i|$, $i=2,\ldots,m$ and that $h_1=1$.

It follows by Corollary~\ref{atomleft} that $$|U_1\cap gU_1|\leq
n-1,$$ thus $|U_1|\leq n$. In addition, for $i\geq 2$,  we have
$$|U_1\cap h_i^{-1}U_i|\leq n-1,$$
 thus $|U_i|\leq n-1$. Therefore $|U|\le m(n-1)+1$ and
$$
|U\cap Ug|= |U|-m\leq \frac{n-2}{n-1}|U|+\frac{1}{n-1}
\le \frac{n-1}n\,|U|,
$$
as claimed. \ Q.E.D.

The minimality of the cardinality of atoms directly yields (see
\cite{Ham96} or \cite{HLS98})

\begin{lem}
\label{atomnonunique}
If $U$ is an $n$-atom
for $C\subset G$, $n\geq 1$, in a torsion-free group $G$,
and $|U|>n$, then any element in $UC$ can be represented
in at least two ways as a product of an element of $U$
and an element of $C$.
\end{lem}

We deduce two rough, but useful estimates about
2-atoms which can
be found in \cite{Ham96}  as well.

\begin{lem}
\label{2atomrough}
If $U$ is a $2$-atom for $C\subset G$, $|C|\geq 3$, in a torsion-free 
group $G$, then $|U|\leq |C|-1$.
\end{lem}
{\it Proof: }
We may assume that $|U|>2$ and $1\in C$, and hence $U\subset UC$.
According to Lemma~\ref{atomnonunique}, for any $u\in U$, there are
 $v_u\in U$ and $c_u\in C\backslash 1$ such that $u=v_uc_u$.
 If $c_u=c_w$ for $u\neq w\in U$, then
 $\{v_w,w\}\subset U\cap v_wv_u^{-1}U$, contradicting
Corollary~\ref{atomleft}.
 Therefore $u\mapsto c_u$ is an injective map from $U$
 into $C\backslash 1$. \ Q.E.D.

For any non-empty $C\subset G$, let $C^{-1}=\{g^{-1}:\,g\in C\}$.

\begin{lem}
\label{2atom}
If $U$ is a $2$-atom
for $C\subset G$, $|C|\geq 3$, in a torsion-free group $G$,
and $|UC|\leq |U|+|C|+k$, then
$|U|\leq k+3$.
\end{lem}
{\it Proof: }
 Let $V$ be a $2$-atom for $U^{-1}$ with $1\in V$, thus
$$
|VU^{-1}|-|V|-|U^{-1}|\leq |C^{-1}U^{-1}|-|C^{-1}|-|U^{-1}|=
|UC|- |U|-|C|\leq k.
$$
If $V=\{1,g\}$ with $g\neq 1$, then Lemma~\ref{atomright} yields
$$
2|U|-1\leq |UV^{-1}|=
|VU^{-1}|\leq k+2+|U|,
$$
which in turn implies $|U|\leq k+3$. If $|V|\geq 3$, then
Lemma~\ref{atomright} and Lemma~\ref{2atomrough} yield
$$
|U|+(|U|-1)+(|U|-2)\leq |UV^{-1}|=|VU^{-1}|\leq k+|V|+|U|\leq k+2|U|-1,
$$
which in turn implies $|U|\leq k+2$. \ Q.E.D.

Now we extend Lemma~\ref{2atom} to $n$-atoms, which extension is the only
novel result of this section.

\begin{prop}
\label{natom} If $U$ is
an $n$-atom for $C\subset G$, $|C|\geq 3$ and
$n\geq 3$, in a torsion-free group $G$, and $|UC|\leq |U|+|C|+k$,
then $|U|\leq n(2k+3)$.
\end{prop}
{\it Proof: } Let $V$ be a $2$-atom for $U^{-1}$,
hence
$$
|VU^{-1}|-|V|-|U^{-1}|\leq |C^{-1}U^{-1}|-|C^{-1}|-|U^{-1}|\leq k.
$$
It follows by Lemma~\ref{2atom} that $|V|\le k+3$. Moreover, by
Lemma~\ref{atomright}, we have $|UV^{-1}|\ge 2|U|-\frac{n-1}{n}|U|$.
Hence,
$$
\mbox{$\frac{n+1}n\,$}|U|\leq |UV^{-1}|=|VU^{-1}|\leq
|U|+|V|+k\leq |U|+2k+3,
$$
thus  $|U|\leq n(2k+3)$. \ Q.E.D.

All these statements about atoms would readily follow from
the following conjecture of Y.O. Hamidoune \cite{Ham09}.

\begin{conj}
\label{atom}
Any $n$-atom in a torsion-free group has
cardinality $n$.
\end{conj}

We recall that a  group $G$ has the unique product property if for
any finite non-empty sets $A,B\subset G$, there is
a $g\in AB$ that
can be represented in a unique way in the form $ab$
with $a\in A$ and $b\in B$. In this case $G$ is torsion-free. It follows by
Lemma~\ref{atomnonunique} that unique product groups satisfy
Conjecture~\ref{atom}.

\section{Small product sets}
\label{sec:proof}

The proof of Theorem \ref{thm:main} together with an estimate of
$c(k)$ will follow from the following Lemma and the
estimations on the size of atoms
in the previous section.

\begin{lem}
\label{uvk} Let $G$ be a torsion-free group. Suppose that $A\subset
G$   with $|A|=3$ is not contained in a left coset of any cyclic
subgroup of $G$.  For $d\geq 3$ and  any finite  set $B\subset G$ of
cardinality
greater than $4d^3$, we have
$$
|AB|> |B|+d.
$$
\end{lem}
{\it Proof: }
We suppose that $|AB|\leq  |B|+d$, and seek a contradiction. We may
assume that $A=\{1,u,v\}$, where $\langle u,v\rangle$ is not cyclic.

 For $g\in G$, we write $B_g=B\setminus 
g^{-1}B=\{x\in B\mid gx\notin B\}$. Since
$|B\cup uB|\le |B|+d$, we see that $|B_u|\le d$. Similarly,
$|B_{u^{-1}}|\le d$, as $|u^{-1}B\cup B|=|B\cup uB|\le|B|+d$, and,
of course, $|B_v|, |B_{v^{-1}}| \le d$ also hold. Since $B$ is
finite, for any $x\in B$ the coset $\langle u \rangle x$ must
contain an element of $B_u$, hence the elements of $B$ belong to at
most $d$ cosets of $\langle u \rangle$, and similarly for $\langle v
\rangle$. Therefore there exists an $x_0 \in B$ such that
$$
|B \cap \langle u \rangle x_0 \cap \langle v \rangle x_0| \ge
|B|/d^2 > 4d.
$$
In order to simplify notation, by replacing $B$ with $Bx_0^{-1}$, we
may assume without loss of generality that $x_0=1$. Let $Z=\langle u
\rangle \cap \langle v \rangle$, and $B_0=B\cap Z$. We have
$|B_0|>4d$. Elements of $Z$ are powers of both $u$ and $v$, hence
$Z$ is contained in the center of $H=\langle u, v\rangle$. As
$Z\neq\{1\}$ and $A$ does not generate a cyclic group, we deduce
that $u$ and $v$ do not commute.

We are going to show that $BZ \supseteq H$. Take an element $g\in H$, and let
us choose a word of shortest length $a_na_{n-1}\cdots a_2a_1$, where
each $a_i$ is one of $u$, $u^{-1}$, $v$, $v^{-1}$, in the coset
$gZ$. Then the cosets $Z$, $a_1Z$, $a_2a_1Z$, \ldots,
$a_na_{n-1}\cdots a_1Z$ are pairwise
disjoint. To any $x\in B_0$, we assign the sequence
$S_x=\{a_ia_{i-1}\cdots a_1x\}_{i=0,1,\ldots,n}$, which sequences  
are pairwise disjoint
as $x$ runs through $B_0$. If $a_na_{n-1}\cdots a_1x\not\in B$,
then there is a smallest $i\in\{1,\ldots,n\}$ such that $a_i\cdots
a_1x\not\in B$. It follows that $S_x$ has an element in
$B_{a_i}\subset B_u\cup B_{u^{-1}}\cup B_v\cup B_{v^{-1}}$,
namely, $x$ if $i=1$, and $a_{i-1}\cdots a_1x$ if $i\geq 2$. Since
$|B_0|>4d\ge |B_u|+| B_{u^{-1}}|+|B_v|+|B_{v^{-1}}|$, and $S_x\cap
S_y=\emptyset$ for $x\neq y$ in $B_0$, there exists an $x\in B_0$
such that $a_na_{n-1}\cdots a_1x\in B$. We conclude that $g\in
a_na_{n-1}\cdots a_1x Z\subset BZ$.

Now the index of the central subgroup $Z$ in $H$ is finite (bounded
by $|B|$), so the center has finite index in $H$. According to a
classical theorem of   Schur (see, e.g., Robinson \cite[Theorem
10.1.4]{Rob96}), this implies that the commutator subgroup of $H$ is
finite. If the commutator subgroup is $1$, then $H$ is abelian, and
if the commutator subgroup is non-trivial, then we have some torsion
elements. In any case, we have contradicted the assumptions on $A$
and $G$, and hence proved the lemma. \ Q.E.D.

\noindent {\it Proof of Theorem~\ref{thm:main}: } Without loss of
generality we may assume that $1\in A$. Then $\langle A \rangle$
is not cyclic by our assumption. Let $B\subset G$ be a finite set
with $|B|> 32(k+3)^6$.
 If $B$ is contained in some right coset of a cyclic subgroup $H$,
then $A$ intersects at least two left cosets of $H$. Let $A_1$
be one of these intersections. Then
using \eqref{eq:kemp} we get
$$
|AB|=|A_1B|+|(A\backslash A_1)B|\geq |A|+2|B|-2>|A|+|B|+k.
$$
Therefore we may assume that $B$ is not contained in a right coset
of any cyclic subgroup.

If $|A|\leq k+3$, then let $A_0=A$, and if $|A|> k+3$, then let
$A_0$ be a $(k+3)$--atom
for $B$ with $1\in A_0$.
By definition, $|AB|\ge |A|+|A_0B|-|A_0|$.
Proposition~\ref{natom} gives that either $|A_0B|-|A_0|>|B|+k$, or
$|A_0|\leq (k+3)(2k+3)$.
If $\langle A_0\rangle$ is not cyclic, then choose $u,v\in A_0\setminus 1$ such
that $\langle u,v\rangle$ is not cyclic. For
$\widetilde{A}_0=\{1,u,v\}\subset A_0$, we have, by Lemma~\ref{uvk},
that
$$
|AB|\geq |A|+|A_0B|-|A_0|\geq
|A|+|\widetilde{A}_0B|-|A_0|
>|A|+(|B|+2(k+3)^2)-|A_0|> |A|+|B|+k.
$$

Finally if $\langle A_0 \rangle$ is cyclic, then
$A_0\ne A$, so $|A_0|\ge k+3$, and $B$ intersects at least two right cosets of
$\langle A_0 \rangle$.
Let $B_1$
be one of these intersections. We have
by \eqref{eq:kemp}
$$
|AB|\geq |A|+|A_0B|-|A_0|=|A|+
|A_0B_1|+|A_0(B\backslash B_1)|-|A_0|
\geq |A|+|B|+|A_0|-2>|A|+|B|+k,
$$
completing the argument. Q.E.D.\\

If $G$ is a unique product group, then
the argument above, just
using Conjecture~\ref{atom} in place
of  Proposition~\ref{natom}, leads to

 \begin{lem}
\label{unique} Let $G$ be a unique product group, $A, B \subset G$
finite subsets, and $k\ge 1$. Suppose that $A$ is not contained in
a left coset of any cyclic  subgroup, and  $|B|> 4(2k+3)^3$, then
  $$
  |AB|> |A|+|B|+k.
  $$
  \end{lem}

\noindent{\bf Acknowledgement: } We are grateful for the help of
Mikl\'os Ab\'ert, Warren Dicks, G\'abor Elek and
Imre Ruzsa in the preparation of this manuscript. We particularly
thank Yahya O. Hamidoune for fruitful discussions on the problem
addressed in this paper, and Peter A. Linnell for
providing in depth information on unique
product groups.

\noindent K\'aroly J. B\"or\"oczky, {\it carlos@renyi.hu}\\
Alfr\'ed R\'enyi Institute of Mathematics, and\\
Universitat Polit\`ecnica de Catalunya, Barcelona Tech, and\\
 Department of Geometry, Roland E\"otv\"os University \\

\noindent P\'eter P. P\'alfy, {\it ppp@renyi.hu}\\
Alfr\'ed R\'enyi Institute of Mathematics\\

\noindent Oriol Serra, {\it oserra@ma4.upc.edu}\\
Universitat Polit\`ecnica de Catalunya, Barcelona Tech

\end{document}